\documentclass[12pt]{article}

\usepackage{amsmath}
\usepackage{amsfonts}
\usepackage{amssymb}
\usepackage{graphicx, color}
\usepackage[colorlinks]{hyperref}
\usepackage{epsfig}
\usepackage{epstopdf} 
\usepackage{enumitem}
\usepackage{mathtools}
\usepackage{comment}

\definecolor{darkgreen}{rgb}{0,0.55,0}

\newtheorem{proposition}{Proposition}[section]
\newtheorem{theorem}{Theorem}[section]
\newtheorem{lemma}[theorem]{Lemma}
\newtheorem{corollary}[theorem]{Corollary}
\newtheorem{remark}[theorem]{Remark}

\newtheorem{definition}{Definition}

\def\phi{{\varphi}}

\DeclareSymbolFont{AMSb}{U}{msb}{m}{n}
\DeclareMathSymbol{\N}{\mathbin}{AMSb}{"4E}
\DeclareMathSymbol{\Z}{\mathbin}{AMSb}{"5A}
\DeclareMathSymbol{\R}{\mathbin}{AMSb}{"52}
\DeclareMathSymbol{\Q}{\mathbin}{AMSb}{"51}
\DeclareMathSymbol{\I}{\mathbin}{AMSb}{"49}
\DeclareMathSymbol{\C}{\mathbin}{AMSb}{"43}

\DeclareMathOperator*{\esssup}{ess\,sup}

\newcommand{\e}{\varepsilon}

\newcommand{\calA}{{\mathcal A}}
\newcommand{\calH}{{\mathcal H}}

\begin{document}
\title{Existence and structure of P-area minimizing surfaces in the Heisenberg group}

\author{ {Amir Moradifam\footnote{Department of Mathematics, University of California, Riverside, California, USA. E-mail: amirm@ucr.edu. Amir Moradifam is supported by NSF grants DMS-1715850 and DMS-1953620.}}\qquad
Alexander Rowell \footnote{Department of Mathematics, University of California, Riverside, California, USA. E-mail: arowe004@ucr.edu.   }
}

\date{\today}
\smallbreak \maketitle
\begin{abstract}
We study existence and structure of $P-$area minimizing surfaces in the Heisenberg group under Dirichlet and Neumann boundary conditions. We show that there exists an underlying vector field $N$ that characterized existence and structure of  $P$-area minimizing surfaces. This vector field exists even if there is no $P$-area minimizing surface satisfying the prescribed boundary conditions. We prove that if $\partial \Omega$ satisfies a so called Barrier condition, it is sufficient to guarantee existence of such surfaces. Our approach is completely different from previous methods in the literature and makes major progress in understanding existence of $P$-area minimizing surfaces.

\end{abstract}

\section{Introduction and Statement of Results}

The $p$-minimal surfaces (also called H-minimal or $X$-minimal surfaces \cite{FSS, GN, P}) have been studied by many authors.  Numerous interesting results have been presented about existence, uniqueness, and regularity of $p$-minimal surfaces \cite{B, CHP, CH, CHMY, CHY, FSS, GN,P}.

Let $\Omega$ be a bounded region in $\R^{2n}$, and $X=(x_1,x'_1,x_2,x'_2, \dots, x_n,x'_n) \in \Omega$. Let $u: \R^{2n}\rightarrow \R$, and consider the graph $(X,u(X))$ in the Heisenberg group of dimension $2n+1$ with prescribed $p$-mean curvature $H(X)$. The $u$ satisfied the equation
\begin{equation}\label{mainPDE}
\nabla \cdot \left( \frac{\nabla u- X^*}{|\nabla u -X^*|} \right)=H,
\end{equation}
where $X^*=(x'_1,-x_1,x'_2,-x_2, \dots, x'_n, -x_n)$ (see Section 2 in \cite{CHY} for a geometric interpretation). The equation \ref{mainPDE} is the Euler-Lagrange equation to the energy functional 
\begin{equation}\label{mainFunctional}
\mathbb{E}(u)=\int_{\Omega} \left(|\nabla u -X^*|+Hu \right)dx_1 \wedge dx'_1 \wedge \dots \wedge dx_n \wedge dx'_n.
\end{equation}
One of the main challenges in studying the equation \eqref{mainPDE} is to deal with the singular set of solutions, i.e. 
\[\{X \in \Omega: |\nabla u(X) -X^*|=0\}.\]
On the other hand, since the energy functional $\mathbb{E}$ is not strictly convex, analysis of existence and uniqueness of minimizers is also a highly non-trivial problem. 

In \cite{B} the author studied the size of the singular set of solutions, and showed the existence of solutions with large singular sets. In \cite{CHY}, the authors proved existence of minimizers of \eqref{mainFunctional} in the special case $H\equiv 0$, and under the assumption that $\Omega$ is a $p$-convex domain (see Theorem A in \cite{CHY}).  They also proved interesting uniqueness and comparison results for minimizers of \eqref{mainFunctional} (Theorem B and C in \cite{CHY}).  In \cite{PSTV}, the authors proved existence and uniqueness of minimizers of $\mathbb{E}$ for the case when $H\equiv 0$ in $\Omega$. In \cite{CH}, the authors studied uniqueness of minimizers of the functional $\mathbb{E}$, and proved several interesting results. 

In this paper, we study existence and structure of minimizers of the energy functional $\mathbb{E}$ from a different point of view, using the Rockafellar-Fenchel duality.  We prove various existence results that are new, even for the case $a\equiv 1$.  Consider the following weighted form of the functional $\eqref{mainFunctional}$ 
\begin{equation}\label{mainFunctionalWeighted}
\mathcal{F}(u)=\int_{\Omega} \left(a |\nabla u -X^*|+Hu \right)dx_1 \wedge dx'_1 \wedge \dots \wedge dx_n \wedge dx'_n,
\end{equation}
where $a \in L^{\infty}(\Omega)$ is a positive function. Minimizers of this functional will satisfy the Euler-Lagrange equation 
\begin{equation}\label{EL1}
    \nabla \cdot \left(a \frac{\nabla u-X^*}{|\nabla u-X^*|} \right)=H,
\end{equation}
which could be viewed as the $p$-mean curvature of the function $(X, u(X))$, with respect to the metric $g=a^{\frac{2}{n-1}}dx$, which is conformal to the Euclidean metric. Our approach is completely different from the previous ones in the literature and provides major progress in understanding the existence of P-area minimizing surfaces.

The equation \eqref{EL1} with $X\equiv  0$ and $H\equiv 0$ have been extensively studied by many authors, including the first author, see \cite{HMN, JMN, Mo, Mo1, MNT, MNTa_SIAM, NTT07, NTT08, NTT10, NTT11,  sternberg_ziemer92, sternberg_ziemer, sternbergZiemer93, ST}.

This paper is organized as follows. In Section 2, we prove existence results under the Neumann boundary condition. In Section 3, we study existence of P-area minimizing surfaces with Dirichlet boundary condition.  Finally, in Section 4 we provide existence of P-area minimizing surfaces under a Barrier condition on the boundary $\partial \Omega$.

\section{Existence of P-area minimizing surfaces with Neumann boundary condition} \label{dual}
Let $\Omega$ be a bounded open region in $\R^n$, $a \in L^{\infty}(\Omega)$ be a positive function, $f\in L^1(\partial \Omega)$, and consider the minimization problem 
\begin{equation}\label{functionalMain0}
	\inf _{u\in \mathring{BV} (\Omega)} I(u):=\int_{\Omega} a \left| D u + F \right| +Hu,
\end{equation} 
where 
\[\mathring{BV}(\Omega)=\{u\in BV(\Omega): \int_{\Omega}u=0\}.\]

In order to study the minimizers of the least gradient problem \eqref{functionalMain} we first analyze the dual of this problem using Rockefeller-Fenchel duality. Define $E(b):L^2(\Omega) \rightarrow \R $ and $G(u): \mathring{H}^1(\Omega) \rightarrow \R$ as follows
$$E(b)=\int_{\Omega} a \left| b+F \right|  \hspace{0.5cm} \text{and} \hspace{0.5cm} G(u)=\int_{\Omega}Hu,$$
where $ \mathring{H}(\Omega)=\{u\in H^1(\Omega): \int_{\Omega}u=0\}$.
Then \eqref{functionalMain} can be rewritten as 
\begin{equation}\label{PrimalPronlem0}
	(P) \hspace{0.5cm} \inf_{u \in  \mathring{H}^1(\Omega)} \{  E(\nabla u)+ G(u)\}.  
\end{equation}

By Rockafellar-Fenchel duality \cite{ET}, the dual problem associated to \eqref{PrimalPronlem0} is 
\begin{equation}\label{dualProb0}
(D) \hspace{0.5cm} -\min_{b \in (L^{2}(\Omega))^n} \{ E^{*}(b)+  G^{*}(-\nabla^* b) \} =\max_{b \in (L^{2}(\Omega))^n} \{-E^{*}(b)-G^{*}(-\nabla^{*}b)\},
\end{equation}
where $E^*$ and $G^*$ are the convex conjugates of the convex functions $E$ and $G$, and $\nabla ^*$ is the adjoint of the gradient operator $\nabla: \mathring{H}^1(\Omega) \rightarrow L^2(\Omega)$. Let us first compute $G^*(-\nabla^* b)$. 

\begin{align*}
	G^*(-\nabla^* b) &=\sup_{u \in H^1(\Omega)} \left\{ \langle u, -\nabla^* b \rangle_{\mathring{H}^1 (\Omega) \times (H^1(\Omega))^*} -G(u) \right\}\\ 
	&=\sup_{u \in \mathring{H}^1(\Omega)} \left\{ \langle u, -\nabla^* b \rangle_{H^1 (\Omega) \times (H^1(\Omega))^*} -\int_{\Omega}Hu \right\} \\
	&=\sup_{u \in \mathring{H}^1(\Omega)} \left\{ -\int_{\Omega} \nabla u \cdot b -\int_{\Omega}Hu \right\}.
\end{align*}  
Since $cu \in  \mathring{H}^1$ for any $u \in \mathring{H}^1$ and any $c\in \R$, 

\begin{equation}
G^*(-\nabla^*b)=
	\begin{cases}
	\infty & \text{ if } u \not \in \mathcal{D}_0(\Omega)\\
	0 & \text{ if } u \in \mathcal{D}_0,
	
	\end{cases}
\end{equation}
where 
\begin{equation}
    \mathcal{D}_0:=\left \{b\in (L^2(\Omega))^n: \int_{\Omega} \nabla u \cdot b+Hu =0, \ \ \hbox{for all}\ \ u\in  \mathring{H}^1(\Omega)\right \}.
\end{equation}
On the other hand, it follows from Lemma 2.1 in \cite{MNT} that 

\begin{equation}
E^*(b)=
	\begin{cases}
	-\langle F,b \rangle & \text{ if } |b|\leq a \ \ \hbox{a.e.} \ \ \hbox{in}\ \ \Omega\\
	\infty & \text{ otherwise }.
	\end{cases}
\end{equation}

Thus, the dual problem $(D)$ can be written as 
\begin{equation}
   (D) \hspace{0.5cm} \sup \{\langle F,b \rangle: b\in \mathcal{D}_0 \ \ \hbox{and}\ \ |b| \leq a \ \ \hbox{a.e. in}\ \ \Omega \}. 
\end{equation}

Let $\nu_\Omega$ denote the outer unit normal vector to $\partial\Omega$. Then for every $b \in(L^{\infty}(\Omega))^n$ with $\nabla\cdot b \in L^n(\Omega)$ there exists a unique function $[b,\nu_\Omega]\in L^{\infty}_{\mathcal{H}^{n-1}}(\partial\Omega)$ such that 
\begin{equation}\label{trace}
\int_{\partial\Omega}[b,\nu_\Omega]u\,d \mathcal{H}^{n-1}=\int_\Omega u\nabla\cdot b dx+\int_\Omega b \cdot D udx,
\quad u\in C^1(\bar\Omega).
\end{equation}
Moreover, for $u\in BV(\Omega)$ and $b \in(L^{\infty}(\Omega))^n$ with $\nabla \cdot b \in L^n(\Omega)$, the linear functional $u\mapsto(b \cdot Du)$ gives rise to a Radon measure on $\Omega$, and \eqref{trace} is valid for every $u\in BV(\Omega)$ (see \cite{Al, An} for a proof). The following lemma is an immediate consequence of \eqref{trace}.

\begin{lemma} \label{LemmaDual}
Let $b \in (L^{\infty}(\Omega))^n \cap \mathcal{D}_0.$ Then 
\[\nabla \cdot b= H-\int_{\Omega}H dx \ \ \hbox{a.e. in}\ \ \Omega,\]
and
\[[b,\nu_{\Omega}]=0 \ \ \mathcal{H}^{n-1}-a.e. \ \ \hbox{ on } \ \ \partial \Omega.\]
\end{lemma}

Indeed, it follows from the above lemma that for any solution $N$ of the dual problem (D) $\nabla \cdot N=H $ a.e. in $\Omega$, and $N$ is orthogonal to the unit normal vector on $\partial \Omega$ in a weak sense. We are now ready to present the main result of this section.

\begin{theorem} \label{Structure}
Let $\Omega$ be a bounded domain in $\R^n$, $F,H \in L^{2}(\Omega)$, and $a \in L^{2}(\Omega)$ be a positive function. Then the duality gap is zero and the dual problem $(D)$ has a solution, i.e. there exists a vector field $N \in \mathcal{D}_0$ with $|N|\leq a$, $|Du+F|-a.e.$ in $\Omega$, such that 
\begin{equation}\label{dualityGap}
\inf _{u\in \mathring{H}^1(\Omega)} \int_{\Omega} \left( a \left| D u + F \right| +Hu \right) dx= \langle F, N\rangle
\end{equation}
Moreover 
\begin{equation}\label{directionParallel}
\frac{Du+F}{|Du+F|}= N, \ \ \ \ |Du+F|-a.e. \ \ \hbox{in}\ \ \Omega,
\end{equation}
for any minimizer $u$ of \eqref{PrimalPronlem0}.  
\end{theorem}
{\bf Proof.} It is easily verified that  $I(v)=\int_{\Omega} a|Dv|$ is convex, and $J: (L^2(\Omega))^n\rightarrow \R$ with $J(p)=\int_{\Omega} a|p| dx$ is continuous at $p=0$. Hence, it follows from Theorem III.4.1 in \cite{ET} that the duality gap is zero and the dual problem $(D)$ has a solution $N$, and consequently \eqref{dualityGap} holds.  

Now let $u\in \mathcal{H}$ be a minimizer of \eqref{PrimalPronlem0}. Then 
\begin{eqnarray*}
\langle F, N\rangle &=& \int_{\Omega} a \left| D u + F \right| + \int_{\Omega} Hu\\
&\geq & \int_{\Omega} |N| |Du+F|+\int_{\Omega} Hu \\ 
&\geq &  \int_{\Omega} N \cdot (Du+F)+\int_{\Omega} Hu \\ 
&=& \langle F,N \rangle+\int_{\Omega} N \cdot Du+Hu\\
&=&  \langle F,N \rangle, 
\end{eqnarray*}
since $N \in \mathcal{D}_0$. Therefore, both the inequalities above are equalities, and hence \eqref{directionParallel} holds. \hfill $\Box$

\begin{remark}
The primal problem $(P)$ may not have a minimizer in $H^1$, but the dual problem $(D)$ always has a solution $N \in (L^2(\Omega))^n$. Note also that the functional $I(u)$ is not strictly convex, and it may have multiple minimizers (see \cite{JMN}). Theorem \ref{Structure} asserts that if $u_1$ and $u_2$ are both minimizers of $(P)$, then 
\begin{equation}
\frac{Du_1+F}{|Du_1+F|}(x)=\frac{Du_2+F}{|Du_2+F|}(x)= N(x),
\end{equation}
for a.e. point $x\in \Omega$ where $|Du_1+F|$ and $|Du_2+F|$ do not vanish. 
\end{remark}

Next we show that if the primal problem $(P)$ is bounded below, then it has a solution in $BV(\Omega)$. The proof follows from standard facts about $BV$ functions, and we sketch it out for the sake of completeness. 

\begin{proposition}\label{ExistenceProp} There exists a constant $C$, depending on $\Omega$, such that if 

	\begin{equation}\label{BoundExistence}
		\max_{x \in \overline{\Omega}} |H(x)|<C,
	\end{equation}
	then the primal problem (P) has a minimizer.
\end{proposition}
\textbf{Proof.}  Let $u_n$ be the minimizing sequence for $I(u)$. Then 
$$ \int |\nabla u_n | - \int |F| -\int |H| |u_n| \leq \int |\nabla u_n | - \int F +\int H u_n \leq \int |\nabla u_n + F| + Hu_n < c,$$
for some constant $c$ independent of $n$. Hence 
$$\int |\nabla u_n | \leq C+ \int |H| |u_n| + \int |F|. $$
It follows from the Poincar\'e's inequality that there exists  a constant $C_{\Omega}$, independent of $n$, such that 
$$ \int |\nabla u_n | \leq C+ ||H ||_{L^{\infty}(\Omega)}  C_{\Omega}\int  |\nabla u_n| + \int |F| $$
$$\Rightarrow \left( 1-C_{\Omega}||H ||_{L^{\infty}(\Omega)}  \right) \int |\nabla u_n| \leq C + \int |F|.  $$

$$\int |\nabla u_n | \leq C'= \frac{C+\int |F|}{\left( 1- C_{\Omega}||H ||_{L^{\infty}(\Omega)}  \right)}$$
provided that $1-C_{\Omega}||H ||_{L^{\infty}(\Omega)} >0$ or equivalently 
\[||H ||_{L^{\infty}(\Omega)} \leq C:=\frac{1}{C_{\Omega}}.\]

It follows from standard compactness results for $BV$ functions that $u_n$ has a subsequence, denoted by $u_n$ again, such that $u_n $ converges strongly in $L^1$ to a function $\hat{u}\in BV$, and $Du_n$ converges to $D\hat{u}$ is the sense of measures. Since the functional $I(u)$ is lower semicontinuous, $\hat{u}$ is a solution of the primal problem \eqref{functionalMain0}. \hfill $\Box$

\begin{corollary}
Let $\Omega$ be a bounded domain in $\R^n$, $F \in L^{2}(\Omega)$ and $a \in L^{2}(\Omega)$ be a positive function. There exists a constant $C$ such that if $||H ||_{L^{\infty}(\Omega)} <C$, then the equation 
\[\nabla \cdot \left ( \frac{Du+F}{|Du+F|}\right)=H - \int_{\Omega}H\]
has a solution $u \in BV(\Omega)$, i.e. there exists $N \in \mathcal{D}_0$ such that 
\[\frac{Du+F}{|Du+F|}=N.\]

\end{corollary}

\section{Existence of P-area minimizing surfaces with Dirichlet boundary condition}
In this section we study existence of $p$-area minimizing surfaces with a given Dirichlet boundary condition on the boundary $\partial \Omega$. Let $\Omega$ be a bounded open region in $\R^n$, $a \in L^{\infty}(\Omega)$ be a positive function, $f\in L^1(\partial \Omega)$, and consider minimization problem 
\begin{equation}\label{functionalMain}
	\inf _{u\in BV_f(\Omega)} I(u):=\int_{\Omega} a \left| D u + F \right| +Hu,
\end{equation} 
where 
\[BV_f(\Omega)=\{u\in BV(\Omega): u|_{\partial \Omega}=f\}.\]
The function $f \in L^1(\partial \Omega)$ can be extended to a function in $W^{1,1}(\Omega)$ (denoted by $f$ again), and the weighted least gradient  problem \eqref{functionalMain} can be written as 
\begin{equation*}
	\inf _{u\in BV_0(\Omega)} I(u):=\int_{\Omega} a \left| D u + \tilde{F} \right| +Hu+\int_{\Omega}Hf dx,
\end{equation*} 
where $\tilde{F}=F+\nabla f$, and $\int_{\Omega}Hf dx$ is a constant. Hence the minimization problem \eqref{functionalMain} is equivalent to the least gradient problem 
\begin{equation}\label{functionalMainZero}
	\inf _{u\in BV_0(\Omega)} I(u):=\int_{\Omega} a \left| D u + F \right| +Hu. 
\end{equation} 
It is easy to verify that the minimizers of \eqref{functionalMain} in $BV_0(\Omega)$ satisfy the Euler-Lagrange equation 
\begin{equation}\label{EL}
    \nabla \cdot \left(a \frac{D u+F}{|D u+F|} \right)=H,
\end{equation}
with $u|_{\partial \Omega}=0.$ However, the minimization problems \eqref{functionalMain} and \eqref{functionalMainZero} do not necessarily have minimizers even if they are bounded below. This is in contrast with our results in Section 2 where boundedness of the functional $I(u)$ in \eqref{functionalMain0} from below automatically implies existence of a minimizer. To see this suppose $u_n$ is a minimizing sequence for \eqref{functionalMainZero} that converges in $L^1(\Omega)$ to a function $\hat{u} \in BV(\Omega)$. Then it follows from lower semicontinuity of the functional $I(u)$ that 
\[I(\hat{u})\leq \inf_{u \in BV_0(\Omega)} I(u).\]
However, the trace of $\hat{u}$ on $\partial \Omega$ may not necessarily be equal to zero. This is the main reason for nonexistence of minimizers for  \eqref{functionalMainZero}. Indeed it is well known that \eqref{functionalMainZero} may not have a minimizer, and proving existence of minimizers for \eqref{functionalMainZero} is a challenging problem that we aim to tackle in this section.  

Similar to the our approach in Section 2, we first analyze the dual of the minimization problem \eqref{functionalMain0} which will be a crucial tool in our analysis. 

\subsection{The Dual Problem}

As in Section 2, let $E(b): (L^2(\Omega))^n \rightarrow \R $ and $G(u):H^1_0(\Omega) \rightarrow \R$ as 
$$E(b)=\int_{\Omega} a \left| b+F \right|  \hspace{0.5cm} \text{and} \hspace{0.5cm} G(u)=\int_{\Omega}Hu,$$
we can rewrite \eqref{relaxedProblem} as 

\begin{equation}\label{PrimalNew}
	(P') \hspace{0.5cm} \inf_{u \in H^1_0(\Omega) } \{  E(\nabla u)+ G(u)\}.  
\end{equation}

By Rockafellar-Fenchel duality \cite{ET}, the dual problem associated to \eqref{PrimalNew} is 
\begin{equation}\label{dualNew}
(D') \hspace{0.5cm} -\min_{b \in (L^{2}(\Omega))^n} \{ E^{*}(b)+  G^{*}(-\nabla^* b) \} =\sup_{b \in (L^{2}(\Omega))^n} \{-E^{*}(b)-G^{*}(-\nabla^{*}b)\},
\end{equation}
where $E^*$ and $G^*$ are the convex conjugates of the convex functions $E$ and $G$, and $\nabla ^*$ is the adjoint of the gradient operator $\nabla: H^1_0(\Omega)  \rightarrow L^2(\Omega)$. Due to the change in our function space, we update the computation of $G^*(-\nabla^* b)$. 
\begin{align*}
	G^*(-\nabla^* b) &=\sup_{u \in H^1_0(\Omega) } \left\{ \langle u, -\nabla^* b \rangle_{H^1_0(\Omega)  \times (H^1_0(\Omega) )^*} -G(u) \right\}\\ 
	&=\sup_{u \in H^1_0(\Omega) } \left\{ \langle u, -\nabla^* b \rangle_{H^1_0(\Omega)  \times (H^1_0(\Omega) )^*} -\int_{\Omega}Hu \right\} \\
	&=\sup_{u \in H^1_0(\Omega) } \left\{ -\int_{\Omega} \nabla u \cdot b -\int_{\Omega}Hu \right\}.
\end{align*}  
Since $cu \in H^1_0(\Omega) $ for any $u\in H^1_0(\Omega) $ and any $c\in \R$, 

\begin{equation}
G^*(-\nabla^*b)=
	\begin{cases}
	\infty & \text{ if } u \not \in \widetilde{\mathcal{D}}_0(\Omega)\\
	0 & \text{ if } u \in \widetilde{\mathcal{D}}_0,
	
	\end{cases}
\end{equation}
where 
\begin{equation}
    \widetilde{\mathcal{D}}_0:=\left \{b\in (L^2(\Omega))^n: \int_{\Omega} \nabla u \cdot b+Hu =0, \ \ \hbox{for all}\ \ u\in H^1_0(\Omega)   \right \}\subseteq \mathcal{D}_0.
\end{equation}
On the other hand, it follows from Lemma 2.1 in \cite{MNT} that 

\begin{equation}
E^*(b)=
	\begin{cases}
	-\langle F,b \rangle & \text{ if } |b|\leq a \ \ \hbox{a.e.} \ \ \hbox{in}\ \ \Omega\\
	\infty & \text{ otherwise }.
	\end{cases}
\end{equation}

Thus the dual problem $(D')$ can be written as 
\begin{equation}
   (D') \hspace{0.5cm} \sup \{\langle F,b \rangle: b\in \widetilde{\mathcal{D}}_0 \ \ \hbox{and}\ \ |b| \leq a \ \ \hbox{a.e. in}\ \ \Omega \}. 
\end{equation}
It follows from the integration by parts formula \eqref{trace} that $b \in (L^{\infty}(\Omega))^n \cap \widetilde{\mathcal{D}}_0$ if and only if 
\[\nabla \cdot b= H \ \ \hbox{a.e. in}\ \ \Omega.\]
We are now ready to prove the following theorem.

\begin{theorem} \label{Structure2}
Let $\Omega$ be a bounded domain in $\R^n$, $F,H \in L^{2}(\Omega)$, $a \in L^{2}(\Omega)$ be a positive function, and assume $(P')$ is bounded below. Then the duality gap is zero and the dual problem $(D')$ has a solution, i.e. there exists a vector field $N \in \widetilde{\mathcal{D}}_0$ with $|N|\leq a$, $|Du+F|-a.e.$ in $\Omega$, such that 
\begin{equation}\label{dualityGap2}
\inf _{u\in H^1_0(\Omega) } \int_{\Omega} \left( a \left| D u + F \right| +Hu \right) dx= \langle F, N\rangle
\end{equation}
Moreover 
\begin{equation}\label{directionParallel2}
\frac{Du+F}{|Du+F|}= N, \ \ \ \ |Du+F|-a.e. \ \ \hbox{in}\ \ \Omega,
\end{equation}
for any minimizer $u$ of \eqref{PrimalNew}.  
\end{theorem}
{\bf Proof.} It is easily verified that  $I(v)=\int_{\Omega} a|Dv|$ is convex, and $J: (L^2(\Omega))^n\rightarrow \R$ with $J(p)=\int_{\Omega} a|p| dx$ is continuous at $p=0$. Hence, it follows from Theorem III.4.1 in \cite{ET} that the duality gap is zero and the dual problem $(D)$ has a solution $N$, and consequently \eqref{dualityGap2} holds.  

Now let $u\in A_0$ be a minimizer of \eqref{PrimalNew}. Since $N \in \widetilde{\mathcal{D}}_0$, we have  
\begin{eqnarray*}
\langle F, N\rangle &=& \int_{\Omega} a \left| D u + F \right| + \int_{\Omega} Hu\\
&\geq & \int_{\Omega} |N| |Du+F|+\int_{\Omega} Hu \\ 
&\geq &  \int_{\Omega} N \cdot (Du+F)+\int_{\Omega} Hu \\ 
&=& \langle F,N \rangle+\int_{\Omega} N \cdot Du+Hu\\
&=&  \langle F,N \rangle. 
\end{eqnarray*}
Therefore, both the inequalities above are equalities, and \eqref{directionParallel2} holds. \hfill $\Box$

\begin{remark}
Note that the primal problem $(P')$ may not have a minimizer in $H^1_0$, but the dual problem $(D')$ always has a solution $N \in (L^2(\Omega))^n$. Note also that the functional $I(u)$ is not strictly convex, and it may have multiple minimizers (see \cite{JMN}). Theorem \ref{Structure} asserts that if $u_1$ and $u_2$ are both minimizers of $(P)$, then 
\begin{equation}
\frac{Du_1+F}{|Du_1+F|}(x)=\frac{Du_2+F}{|Du_2+F|}(x)= N(x),
\end{equation}
for a.e. point $x\in \Omega$ where $|Du_1+F|$ and $|Du_2+F|$ do not vanish. 
\end{remark}

\subsection{The relaxed problem }
Here we study existence of minimizier for the relaxed least gradient problem 
\begin{equation}\label{relaxedProblem}
\inf_{u\in A_0}I(u)= \inf_{u\in A_0} \int_{\Omega}(a|Du+F|+Hu)dx  +\int_{\partial \Omega} a|u| ds,
\end{equation}
where 
$$A_0 := \left \{  u \in H^1 (\R^n) : u=0 \text{ in } \Omega^c   \right \}. $$
Unlike the problem \eqref{functionalMainZero}, any minimizing sequence for \eqref{relaxedProblem} converges to a minimizer in $A_0$.  Indeed the following proposition holds. 

\begin{proposition} \label{PropLast} There exists a constant $C$, depending on $\Omega$, such that if 

	\begin{equation}
		\max_{x \in \overline{\Omega}} |H(x)|<C,
	\end{equation}
	then the primal problem \eqref{functionalMainZero} has a minimizer in $A_0$.
\end{proposition}
{\bf Proof.} The proof follows from an argument similar to the one used in the proof of Proposition \ref{ExistenceProp}, and the observation that if $u_n \in A_0$ converges to $\hat{u}$ in $L^1(\Omega)$, then $\hat{u}\in A_0$. \hfill $\Box$

The next theorem characterizes the relationship between these two problems and sheds light on the challenging problem of existence of minimizers for \eqref{functionalMainZero}. 

\begin{theorem}\label{twoPrimalProblems}
	Let $\Omega \subset \R ^n$ be a bounded open set with Lipschitz boundary, $F \in (L^2(\Omega))^n$, and $H\in L^2(\Omega)$.  If the minimization problem \eqref{functionalMainZero} is bounded below, then

	\begin{equation}\label{MinimizerA0}
		\min_{u \in A_0} \left( \int_{\Omega} (a|D u +F| +Hu)dx +\int_{\partial \Omega} a|u| ds  \right) = \inf_{\substack{u\in BV_0(\Omega) }} \int_{\Omega} a|D u +F| +Hu 
	\end{equation}
Moreover, if $u$ is a minimizer of \eqref{relaxedProblem}, then 
\begin{equation}\label{BoundaryResult1}
u [N, \nu_{\Omega}]=|u| \ \ \mathcal{H}^{n-1}-a.e. \ \ \hbox{on}\ \ \partial \Omega. 
\end{equation}

\end{theorem}
{\bf Proof.}  Since $BV_0(\Omega)$ can be continuously embedded in $A_0$, we have  

$$\min_{u \in A_0} \left( \int_{\Omega} (a|\nabla u +F| +Hu)dx +\int_{\partial \Omega} a|u| ds \right) \leq \inf_{\substack{u\in BV_0(\Omega) }} \int_{\Omega} a|\nabla u +F| +Hu.$$
It follows from Theorem \ref{Structure2} that there exists a vector field $N$ with $|N|\leq a$ a.e. in $\Omega$ and 
\[N=\frac{D u +F}{\left| D u +F \right|}.\]
Now let $u$ be a minimizer of the relaxed problem with $u|_{\partial \Omega}=g|_{\partial \Omega}$, where $g\in W^{1,1}(\Omega)$. Since $u-g \in \widetilde{\mathcal{D}}_0$, we have  

\begin{align*}
	 \min_{u \in A_0} \left( \int_{\Omega} (a|D u +F| +Hu)dx +\int_{\partial \Omega} a|u| ds  \right)&=
    \int_{\Omega} a|\nabla u +F| +Hu +\int_{\partial \Omega} a|u| \\
    & \geq \int_{\Omega} |N||\nabla u +F| +Hu +\int_{\partial \Omega} a|u|\\
    &\geq \int_{\Omega} N (\nabla u +F) +Hu +\int_{\partial \Omega} a|u|\\
    & = \int_{\Omega} N \cdot F + \int_{\Omega} N \cdot \nabla u +Hu +\int_{\partial \Omega} a|u|\\
    & = \langle N, F \rangle + \int_{\Omega} N \cdot \nabla (u-g) + H(u-g) \\
    & + \int_{\Omega} N \cdot \nabla g +Hg +\int_{\partial \Omega} a|g|\\
    &= \langle N, F \rangle + + \int_{\Omega} N \cdot \nabla g +Hg +\int_{\partial \Omega} a|g|\\     
    &= \langle N, F \rangle + \int_{\partial \Omega} g  [N \cdot \nu_{\Omega}] + \int_{\partial \Omega} a|g|\\     
    & \geq \langle N, F \rangle\\
    & = \inf_{\substack{u\in BV_0(\Omega) }} \int_{\Omega} a|D u +F| +Hu. 
\end{align*}
We used integration by parts, and $|N|\leq a$ a.e. in $\Omega$, to obtain the last inequality, and hence \eqref{MinimizerA0} holds. Moreover, all the inequalities in the above computation are equalities. In particular, \eqref{BoundaryResult1} holds.  $\Box$

The following theorem is an immediate consequence of Theorem \ref{Structure2} and Theorem \ref{twoPrimalProblems}.

\begin{theorem}\label{summaryTheorem}
Let $\Omega$ be a bounded domain in $\R^n$, $F,H \in L^{2}(\Omega)$, $a \in L^{2}(\Omega)$ be a positive function, and assume $(P')$ is bounded below. Then there exists a vector field $N \in \widetilde{\mathcal{D}}_0$ with $|N|\leq a$, $|Du+F|-a.e.$ in $\Omega$, such that 
\begin{equation}
\frac{Du+F}{|Du+F|}= N, \ \ \ \ |Du+F|-a.e. \ \ \hbox{in}\ \ \Omega,
\end{equation}
for any minimizer $u$ of \eqref{functionalMainZero}. Moreover, every minimizer of \eqref{functionalMainZero} is a minimizer of \eqref{relaxedProblem}, and  if $u$ is a minimizer of \eqref{relaxedProblem}, then 
\begin{equation}\label{BoundaryResult}
u [N, \nu_{\Omega}]=|u| \ \ \mathcal{H}^{n-1}-a.e. \ \ \hbox{on}\ \ \partial \Omega. 
\end{equation}
In particular, $u=0$ $\mathcal{H}^{n-1}$ a.e. on the set
\[\{x\in \partial \Omega: [N,\nu_{\Omega}]<|N| \}.\]
\end{theorem}

The next theorem follows immediately from Theorem \ref{summaryTheorem}. \\
\begin{theorem}
Let $\Omega$ be a bounded domain in $\R^n$, $F,H \in L^{2}(\Omega)$, $a \in L^{2}(\Omega)$ be a positive function, and assume $(P')$ is bounded below. Let $N$ be the solution of the dual problem guaranteed by Theorem \ref{Structure2} and assume that $[N,\nu_{\Omega}]<|N|$ almost everywhere on $\partial \Omega$. Then the least gradient problem \eqref{functionalMainZero} has a minimizer in $BV_0(\Omega)$. 
\end{theorem}

\section{Existence of minimizers under the Barrier condition}
Let $F\in (L^1(\Omega)^n)$ and $a, H \in L^{\infty} (\Omega)$ with $a >0$ in $\Omega$, and define $\psi:\R^n \times BV_0(\Omega)$ as follows 
\begin{equation}\label {psi}
 \psi(x, u) := a(x) |Du + F\chi_{E_u}| + Hu,
\end{equation}
where $E_u$ is the closure of the support of $u$ in $\Omega$.   

Define the  $\psi$-perimeter of $E$ in $A$, as
\begin{equation*}
    P_\psi(E;A):= \int_{A} a(x) \left| D \chi_E + F \chi_E \right| + H \chi_E. 
\end{equation*}

\begin{definition} \label{AreaMinimizingDef}
\begin{enumerate}
    \item A function $u \in BV(\R^n)$ is $\psi$-total variation minimizing in $\Omega \subset \R^n$ if
    $$
   \int_{\Omega}\psi(x,u) \leq \int_{\Omega}\psi(x,v) \text{  for all } v \in BV(\R^n) \text{ such that } u= v \text{ a.e. in } \Omega^c.
    $$
    \item  A set $E \subset \R ^n$ of finite perimeter is $\psi$-area minimizing in $\Omega$ if
\[P_\psi (E;\Omega) \leq P_\psi (\tilde{E}) \]

  for all $ \tilde{E} \subset \R^n \text{ such that } \tilde{E} \cap \Omega^c = E \cap \Omega ^c \text{ a.e.}$.
  
\end{enumerate}
\end{definition}

We will show that the super level sets of $\psi$-total variation minimizing functions in $\Omega$ are $\psi$-area minimizing in $\Omega$. In order to achieve this, we shall first prove some preliminary lemmas.

\begin{lemma} \label{LSC-Lemma} Let $\chi_{\epsilon,\lambda}$ be defined as in \eqref{chiEpsilon}. Then 
\[ P_{\psi}(E,\Omega) \leq \liminf_{\epsilon \rightarrow 0} \int_{\Omega} a(x)|D \chi_{\epsilon,\lambda}+F\chi_{\chi_{\epsilon, \lambda}}|+H \chi_{\epsilon, \lambda}. \]
\end{lemma}
{\bf Proof.} We have 
\begin{eqnarray*}
&& \int_{\Omega} a(x)|D \chi_{\epsilon,\lambda}+F\chi_{\chi_{\epsilon, \lambda}}|+H \chi_{\epsilon, \lambda} -\int_{\Omega} a(x)|D \chi_{E}+F\chi_{E}|+H \chi_{E} \\
&=&  \int_{\Omega \cap \{\lambda-\epsilon<u < \lambda+\epsilon\}} a|D \chi_{\epsilon,\lambda}+F\chi_{\chi_{\epsilon, \lambda}}|+H \chi_{\epsilon, \lambda} - a|D \chi_{E}+F\chi_{E}|-H \chi_{E} \\
&\geq &  \int_{\Omega \cap \{\lambda-\epsilon<u < \lambda+\epsilon\}} a|D \chi_{\epsilon,\lambda}|-a|F\chi_{\chi_{\epsilon, \lambda}}|+H \chi_{\epsilon, \lambda} - a(x)|D \chi_{E}|-a|F\chi_{E}|-H \chi_{E} \\
&=&  \int_{\Omega \cap \{\lambda-\epsilon<u < \lambda+\epsilon\}} a|D \chi_{\epsilon,\lambda}|- a(x)|D \chi_{E}|+H \chi_{\epsilon, \lambda} -H \chi_{E} -a|F\chi_{\chi_{\epsilon, \lambda}}|- a|F\chi_{E}|\\
&=& \int_{\Omega } a|D \chi_{\epsilon,\lambda}|- \int_{\Omega}a(x)|D \chi_{E}|+ \int_{\Omega} (H \chi_{\epsilon, \lambda} -H \chi_{E})\\
&& - \int_{\Omega \cap \{\lambda-\epsilon<u < \lambda+\epsilon\}} a|F\chi_{\chi_{\epsilon, \lambda}}|+ a|F\chi_{E}|.\\
\end{eqnarray*}
It is easy to see that the last two integrals converge to zero as $\epsilon \rightarrow 0$. Hence 

\begin{eqnarray*}
&&  \liminf_{\epsilon \rightarrow 0} \int_{\Omega} a(x)|D \chi_{\epsilon,\lambda}+F\chi_{\chi_{\epsilon, \lambda}}|+H \chi_{\epsilon, \lambda}- P_{\psi}(E,\Omega)\\
&= & \liminf_{\epsilon \rightarrow 0}\int_{\Omega} a(x)|D \chi_{\epsilon,\lambda}+F\chi_{\chi_{\epsilon, \lambda}}|+H \chi_{\epsilon, \lambda} -\int_{\Omega} a(x)|D \chi_{E}+F\chi_{E}|+H \chi_{E} \\
&\geq & \liminf_{\epsilon \rightarrow 0} \int_{\Omega } a|D \chi_{\epsilon,\lambda}|- \int_{\Omega}a(x)|D \chi_{E}|\geq 0,\\
\end{eqnarray*}
where we have used the lower semi-continuity of $\int_{\Omega}a|Dv|$ to obtain the last inequality (see \cite{JMN}). The proof is complete. \hfill $\Box$

\vspace{.2cm}

If  $w\in BV(\R^n)$ and $\Omega$ is an open set with Lipschitz boundary, we will write $w^+$ and $w^{-}$  to denote the outer and inner trace of $w$ on $\partial \Omega.$

\begin{lemma}\label{Relaxed-Lemma} Let $\Omega\subset \R^n$ be bounded and open, with Lipschitz boundary.
Given $g\in L^1(\partial \Omega ; \mathcal H^{n-1})$, define 
\[
I_\psi(v ; \Omega, g) := \ \int_{\partial \Omega} a |g - v^- +F_{\chi_v}|
d\mathcal H^{n-1}
+
\int_{\Omega}\psi(x,Dv).
\]
Then $u\in BV(\R^n)$ is $\psi$-total variation minimizing in $\Omega$ 
if and only if $u|_\Omega$ minimizes
$I_\psi ( \, \cdot \, ; \Omega, g)$ for some $g$, and moreover $g = u^+$.
\label{lem:reformulate}\end{lemma}

{\bf Proof:}
First note that if $v\in BV(\R^n)$ then $v^+, v^- \in L^1(\partial \Omega;\mathcal H^{n-1})$, and conversely, for every $g\in L^1(\partial \Omega;\mathcal H^{n-1})$
there exists some $v\in BV(\R^n)$ such that $g = v^+$. Also 
\begin{equation}
\int_{\partial \Omega} \psi(x, Dv) = \int_{\partial \Omega} a |Dv +F_{\chi_v}|
d\mathcal H^{n-1}
=\ \int_{\partial \Omega} a |v^+ - v^- +F_{\chi_v}|
d\mathcal H^{n-1}.
\label{boundary.meas}\end{equation}
To see this, note that $|Dv|$ can only concentrate on a
set of dimension $n-1$ if that set is a subset of the jump set of $v$,
so \eqref{boundary.meas} follows from standard descriptions of
the jump part of $Dv$.

Now if $u, v\in BV(\R^n)$ satisfy $u=v$ a.e. in $\Omega^c$, then  $\int_{\bar \Omega^c} \varphi(x, Du) = \int_{\bar \Omega^c} \varphi(x, Dv)$. In
addition, $u^+ = v^+$, so
using  \eqref{boundary.meas}
we deduce that
\[
\int_{\R^n}\psi(x, Du) - \int_{\R^n} \psi (x,Dv)
 \ = \ I_\varphi(u; \Omega, u^+) - I_\varphi(v ; \Omega, u^+).
\]
The lemma easily follows from the above equality.
\hfill $\Box$
\\

\begin{theorem} \label{SuperLevelThm}
Consider the bounded Lipschitz domain $\Omega \subset \R^n$ and a $\psi$-total variation minimizing function in $\Omega$, $u \in BV(\R^n)$.  Let the super level sets of $u$ to be defined as
\begin{equation}\label{SuperLevelDef}
    E_\lambda := \left \{ x\in \R^n : u(x) \geq \lambda \right \}.
\end{equation}
Then $E_\lambda$ is $\psi$-area minimizing in $\Omega$.

\end{theorem}

{\bf Proof.}  This proof closely mirrors that of Theorem 2.6 in \cite{JMN}.  Consider an arbitrary $\lambda \in \R$, and let $u_1= \max (u- \lambda, 0),  u_2=u-u_1$. For any $g \in BV(\R^n)$ such that $\text{supp}(g) \subset \overline{\Omega},$ we have
\begin{align*}
    \int_{\Omega} a \left| D u_1 + F\chi_{\{u\geq \lambda \}} \right| + H u_1 
    &+ \int_{\Omega} a \left| D u_2 + F\chi_{\{u < \lambda \}} \right| + H u_2 = \int_{\Omega} a \left| D u + F \right| + H u\\
    & \leq \int_{\Omega} a \left| D (u+g) + F \right| + H (u+g) \\
    &=   \int_{\Omega} a \left| D u_1 +D(g\chi_{\{u\geq \lambda\} })+ F\chi_{\{u\geq \lambda \}} \right| + H (u_1+g) \\
    & + \int_{\Omega} a \left| D u_2 +D(g\chi_{\{u< \lambda\} })+ F\chi_{\{u < \lambda \}} \right| + H u_2 \\
     &=   \int_{\Omega} a \left| D u_1 +D(g\chi_{\{u\geq \lambda\} })+ F\chi_{\{u\geq \lambda \}} \right| + H (u_1+g) \\
    & +\int_{\Omega}a|D(g\chi_{\{u< \lambda\} })|+ \int_{\Omega} a \left| D u_2 + F\chi_{\{u < \lambda \}} \right| + H u_2 \\
    & = \int_{\Omega} a \left| D (u_1+g) + F\chi_{\{u\geq \lambda \}} \right| + H (u_1+g)\\
    & + \int_{\Omega} a \left| D u_2 + F\chi_{\{u< \lambda \}} \right| + H u_2.
\end{align*}
Thus
\[ \int_{\Omega} a \left| D u_1 + F\chi_{u_1} \right| + H u_1 \leq \int_{\Omega} a \left| D (u_1+g) + F\chi_{u_1} \right| + H (u_1+g),\]
for all $g \in BV(\R^n)$ with $\text{supp}(g) \subset \overline{\Omega}$. Hence $u_1$ is also $\psi$-total variation minimizing.  By the same process, we can verify that the function defined below is also $\psi$-total variation minimizing,
\begin{equation}\label{chiEpsilon}
    \chi_{\epsilon, \lambda} := \min \left( 1, \frac{1}{\epsilon} u_1  \right) =
    \begin{cases}
        0 & \text{ if } u \leq \lambda, \\
        \frac{1}{\epsilon} (u-\lambda) & \text{ if } \lambda < u \leq \lambda +\epsilon, \\
        1 & \text{ if } u > \lambda + \epsilon.
    \end{cases}
\end{equation}

For a.e. $\lambda \in \R$ the boundary of the super level set $E_\lambda$ is a set of measure zero, that is,
\begin{equation}\label{levelMeasureZero}
     \mathcal{L}^n \left( \{ x \in \Omega : u(x) = \lambda  \} \right) = \mathcal{H}^{n-1} \left( \{ x \in \partial \Omega : u^{\pm}(x) = \lambda  \} \right) =0.  
\end{equation}
 
It follows that 
$$\chi_{\epsilon, \lambda} \rightarrow \chi_\lambda := \chi_{E_\lambda}  \text{ in } L^1_{\text{loc}}(\R^n), \hspace{0.5cm}  \chi_{\epsilon, \lambda}^{\pm} \rightarrow \chi_{\lambda}^{\pm}  \text{ in } L^1(\partial \Omega ; \mathcal{H}^{n-1}), $$
as $\epsilon \rightarrow 0$.

It follows  from Lemma \ref{LSC-Lemma} via quite standard arguments that  
%Let $g := u^+$ on $\partial \Omega$.
%It suffices to show that for any $v\in BV(\Omega)$,
%\[
%I_\varphi(u;\Omega,g) \le I_\varphi(v ; \Omega,g).
%\]
\begin{equation}
P_\psi(\chi_{\lambda}, \Omega) \leq \liminf_{\epsilon \rightarrow 0} P_{\psi}(\chi_{\epsilon,\lambda},\Omega);
\label{wlsc}\end{equation}
and this, with  the $L^1$ convergence of the traces, implies that
\begin{equation}\label{limit1}
I_\varphi(\chi_{\lambda}; \Omega, \chi_{\lambda}^{+}) \le \liminf_{k\to\infty} I_\varphi(\chi_{\epsilon,\lambda}; \Omega, \chi_{\lambda,\epsilon}^{+}).
\end{equation}
Now for any $F\subset \R^n$ such that $\chi_{\lambda}=\chi_F$ a.e. in $\Omega^c$, 
\begin{align*}
I_\varphi(\chi_{\epsilon,\lambda}; \Omega, \chi_{\epsilon,\lambda}^+) 
&\le
I_\varphi(\chi_F; \Omega, \chi_{\epsilon,\lambda}^+)\\
&\le
I_\varphi(\chi_F; \Omega, \chi_\lambda^+) +\int_{\partial \Omega}a|\chi_{\lambda}^+-\chi_{\epsilon,\lambda}^{+}|
 \ d\mathcal H^{n-1}\\
&\le
I_\varphi(\chi_F; \Omega, {\chi_{\lambda}}^+) +C\int_{\partial \Omega}
|\chi_{\lambda}^+-\chi_{\epsilon,\lambda}^{+}|\ d\mathcal H^{n-1}.
\end{align*}
It follows from this, \eqref{limit1},
and $\chi_{\epsilon,\lambda}^+ \rightarrow \chi_{\lambda}^+$ in $L^1(\partial \Omega;\mathcal H^{n-1})$ that

\[I_\varphi(\chi_\lambda;\Omega, \chi_{\lambda}^+) \le 
I_\varphi(\chi_F;\Omega, \chi_\lambda^+),\]
which proves that $E_{\lambda}$ is $\phi$-area minimizing in $\Omega$.

In the case where $\lambda$ does not satisfy (\ref{levelMeasureZero}), we can take an increasing sequence $ \lambda_k  \rightarrow \lambda$ as $k \rightarrow \infty,$ that satisfies (\ref{levelMeasureZero}) for each $k$.  This implies that
$$\chi_{\lambda_k} \rightarrow \chi_\lambda \text{ in } L^1_{\text{loc}}(\R^n), \hspace{0.5cm}  \chi_{\lambda_k}^{\pm} \rightarrow \chi_{\lambda}^{\pm}  \text{ in } L^1(\partial \Omega ; \mathcal{H}^{n-1}). $$
This once again leads to the conclusion that $E_\lambda$ is $\psi$-area minimizing in $\Omega$ in view of Lemma \ref{Relaxed-Lemma}.  \hfill  $\Box$

Now we are ready to present the main existence results of this section. For any measurable set $E$ define

\[E^{(1)}:= \{x\in \R^n :
\lim_{r\to 0} \frac {\calH^n(B(r,x)\cap E)}{\calH^n(B(r))} = 1 \}.\]

\begin{definition} \label{BarrierCondition}
Suppose that $\Omega \subset \R^n$ is a bounded Lipschitz domain. Then $\Omega$ satisfies the barrier condition if for every $x_0 \in \partial \Omega$ and $\epsilon >0$ sufficiently small, if $V$ minimizes $P_\psi ( \cdot ; \R ^n)$ in
\begin{equation}\label{BD-cond}
\{ W \subset  \Omega: W \setminus B(\epsilon, x_0 ) = \Omega \setminus B(\epsilon, x_0 )    \},
\end{equation}
then
$$
\partial V^{(1)} \cap \partial \Omega \cap B(\epsilon, x_0)= \emptyset.
$$
\end{definition}

\begin{lemma} \label{MainLemma}
    Given a bounded Lipschitz domain $\Omega \subset \R^n$ that satisfies the barrier condition from Definition \ref{BarrierCondition}, and suppose $E \subset \R^n$ minimizes $P_\psi(\cdot ;\Omega)$.  Then
    $$
    \left\{ x \in \partial \Omega \cap \partial E^{(1)} : B(\epsilon, x) \cap \partial E^{(1)} \subset \overline{\Omega} \text{ for some } \epsilon > 0 \right\} = \emptyset.
    $$
\end{lemma}
{\bf Proof.} Assume there exists $x_0\in \partial \Omega \cap \partial E^{(1)}$ such that  $B(\epsilon, x_0) \cap \partial E^{(1)}\subset \bar{\Omega}$ for some $\epsilon>0$. Then $\tilde{V}= E\cap \Omega$ is a minimizer of $P_{\psi}(\,\cdot\,; \R^n)$ in (\ref{BD-cond}), and 
\[
x_0\in \partial {\tilde{V} }^{(1)}\cap \partial \Omega \cap B(\epsilon, x_0)\neq  \emptyset. 
\]
This contradicts the  barrier condition and finishes the proof. 
\hfill $\Box$

\medskip

Define
$$ 
BV_f(\Omega) := \left\{ u \in BV(\Omega): \lim_{r \rightarrow 0} \esssup_{y \in \Omega, |x-y|<r} |u(y)-f(y)|=0 \text{ for } x \in \partial \Omega \right\}. 
$$

\begin{theorem} \label{MainThm}
Let $\psi:\R ^n \times \R ^n \rightarrow \R$ be defined as in \eqref{psi}, and $\Omega \subset \R ^n$ be a bounded Lipschitz domain. Suppose $||H||_{L^{\infty}(\overline{\Omega})}$ is small enough such that Proposition \ref{PropLast} holds. If $\Omega$ satisfies the barrier condition with respect to $\psi$, as given in Definition \ref{BarrierCondition}, then for every $f \in C(\partial \Omega)$ the minimization problem \eqref{functionalMain} has a minimizer in $BV_f(\Omega)$. 
\end{theorem}

{\bf Proof.}  Since every $\calH^{n-1}$
integrable function on $\Omega$ is the trace of some
(continuous) function in $BV(\Omega^c)$, without loss of generality we may assume that $f\in BV(\R^n)$.

Define
\[
\calA_f := \{v\in BV(\R^n): \ \ v=f \ \ \hbox{on}\ \ \Omega^c \},
\]
and note that $BV_f(\Omega) \hookrightarrow \calA_f$, 
in the sense that any element $v$ of $BV_f(\Omega)$
is the restriction to $\Omega$ of a unique element of 
$\calA_f$. An argument similar to that of Proposition \ref{PropLast} implies that $\int_{\R^n} \psi(x,v)$ has as a minimizer $u\in \calA_f$.

We next use the barrier condition to show  that $u\in BV_f(\Omega)$.
If not, there exists some $x\in \partial \Omega$
and $\delta>0$ such that
\begin{equation}
 \esssup_{y\in \Omega, |x-y|<r}\big( f(x) - u(y))\ge \delta\qquad
\mbox{ or }  \ \ \esssup_{y\in \Omega, |x-y|<r}\big(u(y) - f(x)) \ge \delta
\label{ess.alt}\end{equation}
for every $r>0$.
Assume that the latter condition holds. 
It follows
from this and the continuity of $f$, 
that $x\in  \partial E^{(1)}$ for $E := E_{ f(x) + \delta/2}$. By Theorem \ref{SuperLevelThm}  $E$ is $\psi$-area minimizing in $\Omega$. 
However, since $f$ is continuous in $\Omega^c$ and
$u\in \calA_f$, it is clear that
$u< f(x) + \delta/2$ in $B(\e, x)\setminus \Omega$ for all sufficiently
small $\e$. This contradicts Lemma \ref{MainLemma}. If the first alternative holds in \eqref{ess.alt}, then we
set $E :=\{ y\in \R^n : u(y) \le f(x) - \delta/2\}$
and reach a similar contradiction. Hence $u\in BV_f(\Omega)$, and in view of Theorem \ref{twoPrimalProblems}, it is $\psi$-total variation minimizing in $BV_f(\Omega)$. \hfill $\Box$

\end{document}